\documentclass[11pt,english]{article}

\RequirePackage[OT1]{fontenc}
\usepackage{graphicx}
\usepackage{epsfig}
\usepackage{epstopdf}
\usepackage{color}
\usepackage{bbm}
\usepackage{amsmath}
\usepackage{tabularx} 
\usepackage{amsthm}
\usepackage{dsfont}
\usepackage{placeins}
\usepackage{amssymb}
\usepackage{float}
\usepackage{mathtools}
\usepackage[square,comma,numbers,sort&compress]{natbib}
\usepackage[numbers]{natbib}
\usepackage{extarrows}
\usepackage{tabularx}
\RequirePackage[colorlinks,citecolor=blue,urlcolor=blue]{hyperref}

\allowdisplaybreaks

\definecolor{gray}{rgb}{0.5,0.5,0.5}
\definecolor{black}{rgb}{0,0,0}
\definecolor{white}{rgb}{1,1,1}
\definecolor{blue}{rgb}{0,0,0.7}

\definecolor{green}{rgb}{0.133,0.545,0.133}

\definecolor{yellow}{rgb}{1,0.549,0}

\definecolor{red}{rgb}{1,0.133,0.133}

\definecolor{purple}{rgb}{0.58,0,0.827}

\definecolor{backgcode}{rgb}{0.97,0.97,0.8}
\definecolor{Brown}{cmyk}{0,0.81,1,0.60}
\definecolor{OliveGreen}{cmyk}{0.64,0,0.95,0.40}
\definecolor{CadetBlue}{cmyk}{0.62,0.57,0.23,0}

\newcommand{\twovec}[2]{\bracks{\begin{array}{c} #1 \\ #2 \end{array}}}

\newcommand{\twobytwomat}[4]{\bracks{\begin{array}{cc} #1 & #2 \\ #3 & #4 \end{array}}}

\newcommand{\beqn}{\vspace{-0.25cm}\begin{eqnarray*}}
\newcommand{\eeqn}{\end{eqnarray*}}
\newcommand{\bneqn}{\vspace{-0.25cm}\begin{eqnarray}}
\newcommand{\eneqn}{\end{eqnarray}}

\newcommand{\beans}{\color{blue} \begin{equation}*  \text{Ans:}~~~}
\newcommand{\eeans}{\end{equation}* \color{black}}

\newcommand{\parens}[1]{\left(#1\right)}

\newcommand{\bracks}[1]{\left[#1\right]}

\newcommand{\abss}[1]{\left|#1\right|}

\newcommand{\expe}[1]{\mathbb{E}\bracks{#1}}

\renewcommand{\sin}[1]{\text{sin}\parens{#1}}

\newtheorem{theorem}{Theorem}

\newtheorem{proposition}[theorem]{Proposition}
\newtheorem{remark}{Remark}

\author{Philip Ernst \footnote{Rice University, Department of Statistics}, Larry Shepp \footnote{Deceased April 23, 2013}, and Abraham Wyner \footnote{The Wharton School of the University of Pennsylvania, Department of Statistics}}
\title{Yule's ``Nonsense Correlation'' Solved!}
\date{\today}

\begin{document}
\maketitle

\begin{abstract}

In this paper, we resolve a longstanding open statistical problem. The problem is to mathematically prove Yule's 1926 empirical finding of ``nonsense correlation'' (\cite{Yule}). We do so by analytically determining the second moment of the empirical correlation coefficient

\beqn
\theta := \frac{\int_0^1W_1(t)W_2(t) dt - \int_0^1W_1(t) dt \int_0^1 W_2(t) dt}{\sqrt{\int_0^1 W^2_1(t) dt - \parens{\int_0^1W_1(t) dt}^2} \sqrt{\int_0^1 W^2_2(t) dt - \parens{\int_0^1W_2(t) dt}^2}},
\eeqn
of two {\em independent} Wiener processes, $W_1,W_2$. Using tools from Fredholm integral equation theory, we successfully calculate the second
moment of $\theta$ to obtain a value for the standard deviation of $\theta$  of
nearly .5. The ``nonsense'' correlation, which we call ``volatile'' correlation, is volatile in the sense that its distribution is heavily dispersed and is frequently large in absolute value. It is induced because each Wiener process is
``self-correlated'' in time. This is because a Wiener process is an integral
of pure noise and thus its values at different time points are correlated.  In addition to providing an explicit formula for the second moment of $\theta$, we offer implicit formulas for higher moments of $\theta$.

\end{abstract}

\section{Introduction}\label{sec:introduction}
A fundamental yet crucial question for practitioners of statistics is the following: given a sequence of pairs of random variables $\{X_k,Y_k \}$ ($k=1,2, \ldots, n$), how can we measure the strength of the dependence of $X$ and $Y$? The classical Pearson correlation coefficient offers a solution that is standard and often powerful. It is also widely used even in situations where little is known about its empirical properties; for example, when the sequence of random variables are not identically distributed or independent.  The Pearson correlation is often calculated between two time series. Practitioners have developed many ``rules of thumb" to help interpret these values (i.e. a correlation greater than .9 is generally understood to be large). Such correlations can be difficult to interpret (\cite{Yule}).  It is well known that a  ``spurious'' correlation will be observed when two time series are themselves dependent on an unobserved third time series. 

However, it is less well known to some practitioners that one may observe ``volatile'' correlation in independent time series. The correlation is volatile in the sense that its distribution is heavily dispersed and is frequently large in absolute value. Yule observed this empirically and, in his 1926 seminal paper (\cite{Yule}), called it ``nonsense'' correlation, asserting that ``we sometimes obtain between quantities varying with time (time-variables) quite high correlations to which we cannot attach any physical significance whatever, although under the ordinary test the correlation would be held to be certainly significant.'' Yet Yule's empirical finding remained ``isolated'' from the literature until 1986 (see \cite{Aldrich}), when the authors of \cite{Hendry} and \cite{phil} confirmed many of the empirical claims of ``spurious regression'' made by the authors of \cite{gran}.  However, despite these advances, mathematical confirmation of ``nonsense'' correlation has remained open. We now turn to closing this matter.

Throughout this work, we will use the word ``volatile'' in lieu of ``nonsense.'' We emphasize that volatile correlation is distinctly different from ``spurious'' correlation, as the latter refers to a third time series, whereas the former does not. Despite Yule's empirical findings, it is often (erroneously) assumed that a large correlation between such pairs must have a cause (see \cite{Mann}) when they definitionally do not. Suppose, for example, that $X_i = S_i$ and $Y_ i= S'_i$ where $S_i$ and $S'_i$ are the partial sums of two independent random walks. The empirical correlation is defined in the usual way as
\begin{equation}
\theta_n = \frac{\sum_{i=1}^n S_i S^\prime_i - \frac{1}{n}(\sum_{i=1}^n S_i)(\sum_{i=1}^n S^\prime_i)}{\sqrt{\sum_{i=1}^n S^2_i - \frac{1}{n}(\sum_{i=1}^n S_i)^2}\sqrt{\sum_{i=1}^n (S^\prime_i)^2 - \frac{1}{n}(\sum_{i=1}^n S^\prime_i)^2}}.
\end{equation}
 Nevertheless, it is sometimes (erroneously) assumed that for large enough $n$ these correlations should be small. Indeed, large values are quite probable.  The histogram of the empirical correlation of independent random walks is widely dispersed over nearly the entire range. This  was recently mentioned by \cite{McShane}, which presents (using $n=149$) a critique of efforts to reconstruct the earth's historical temperatures using correlated time series. We reproduce the histogram here as Figure 1, with $n=10,000$. 

\begin{figure}[H]
\centering
\includegraphics[width=3.5in]{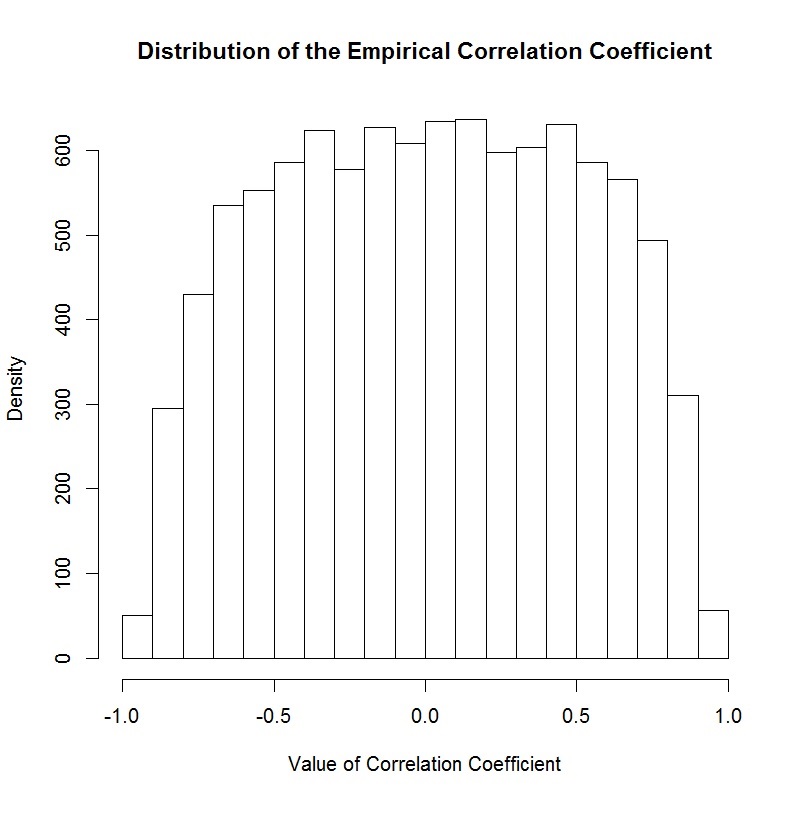} \caption{Simulated empirical correlation of partial sums of two independent random walks with $n=10,000$.} \label{fig:hist} 
\end{figure}
The above histogram reports that the middle 95\% of the observed correlation coefficients fall in the interval [-.83, .83]. The lesson to be learned from Figure \ref{fig:hist} is that correlation is not always a useful statistic; in the case of two independent random walks, the observed correlation coefficient has a very different distribution than that of the nominal $t$-distribution.

\begin{remark}
It should be noted that the statistic which uses the {\em actual random variables}, $X_k, X^\prime_k$
instead of their partial sums, $\sum_{j=1}^k X_j, \sum_{j=1}^k X^\prime_j$ does
not produce volatile correlation. Indeed, as $n \rightarrow \infty$, if the random variables,
$X_k,X^\prime_k, k = 1,\ldots$ are I.I.D. sequences, independent of each other,
with positive finite variances, then

\beqn
\theta^\prime_n = \frac{\frac{1}{n}\sum_{k=1}^n X_k X^\prime_k - \parens{\frac{1}{n}\sum_{k=1}^n X_k}\parens{\frac{1}{n} \sum_{k=1}^n X^\prime_k}}{\parens{\sqrt{\frac{1}{n} \sum_{k=1}^n X_k^2 -(\frac{1}{n} \sum_{k=1}^n X_k)^2}}\parens{\sqrt{\frac{1}{n} \sum_{k=1}^n X_k^2 -(\frac{1}{n} \sum_{k=1}^n X_k)^2}}}
\eeqn
is easily seen to tend to zero (by the law of large numbers). This shows that the volatile correlation is
a consequence of using the partial sums in place of the variables themselves.
The reason that the partial sums are self-correlated (and thereby induce
large correlation) seems related to the arcsine law. The history of Sparre
Andersen's major combinatorial contribution to the proof of the arcsine law
(see \cite{feller}) raises the question of whether a formula 
for discrete sequences of partial sums can be derived by combinatorial methods
employing cyclic permutations. This would be very elegant, and would greatly simplify earlier works of Erd{\"o}s and Kac (see \cite{Kac}), but seems unlikely.
\end{remark}
In lieu of considering the volatile correlation between two independent random walks, we consider the continuous analog, namely that of two independent Wiener processes. Although it would be ideal to find the full analytical distribution of the empirical correlation between two independent Wiener processes, finding the second moment itself suffices as formal confirmation of Yule's ``nonsense'' correlation. But such evidence, until now, has remained elusive.\\
\indent Formally, let $W_i(t), 0\le t \le 1, i = 1,2$ denote two independent Wiener processes on $[0,1]$. We analytically find the second moment of the empirical correlation between $W_1(t)$ and $W_2(t)$, where the empirical correlation is defined as

\begin{equation} \label{corr}
\theta = \frac{\int_0^1W_1(t)W_2(t) dt - \int_0^1W_1(t) dt \int_0^1 W_2(t) dt}{\sqrt{\int_0^1 W^2_1(t) dt - \parens{\int_0^1W_1(t) dt}^2} \sqrt{\int_0^1 W^2_2(t) dt - \parens{\int_0^1W_2(t) dt}^2}}.
\end{equation}
Using tools from integral equation theory, we successfully calculate the second
moment of $\theta$ to obtain a value for the standard deviation of $\theta$  of
nearly .5. This volatile correlation is induced because each Wiener process is
``self-correlated'' in time. This is because a Wiener process is an integral
of pure noise and thus its values at different time points are correlated. The correct intuition for the occurrence of this phenomenon is that a Wiener process {\em self-correlated}, and thus volatile correlation is indeed induced between the independent Wiener processes. \\
\indent Of course, it is rather simple to simulate the distribution of $\theta$.  We do so below using a simple Monte Carlo routine. The result of averaging 10,000 simple Monte Carlo iterations of the first ten moments of $\theta$ can be found below in Table \ref{tab:conditional_prob_sim}.
\small
\begin{table}[ht!]
\centering
\begin{tabular}{ccccccccccccccccccc}
$\expe{\theta^0}$ & $\expe{\theta^1}$ & $\expe{\theta^2}$ & $\expe{\theta^3}$ & $\expe{\theta^4}$& $\expe{\theta^5}$ \\ \hline
1 & -.00116634 & .235057 & -.000524066 & .109276 &-.000283548  \\ \hline
$\expe{\theta^6}$ &$\expe{\theta^7}$ &$\expe{\theta^8}$ &$\expe{\theta^9}$ &$\expe{\theta^{10}}$ \\ \hline
.0609591&-.00016797 &.0378654 &-.000105611 &.0251693 \\ \hline
\end{tabular}
\caption{Moments of $\theta$ obtained from 10,000 Monte Carlo iterations}
\label{tab:conditional_prob_sim}
\end{table}
\normalsize
\noindent The initial validity of the Monte Carlo can be justified by noting that the odd moments should all be zero, since $\theta$ is symmetric. The impact of this problem for practitioners of statistical inference, however, renders it most deserving an analytical solution, and this serves as our paper's primary motivation. \\
\indent The structure of the paper is as follows. In \S \ref{sec:formal_statement_results}, we present results needed for obtaining the distribution of $\theta$. In \S \ref{sec3}, we provide implicit formulas for all moments of
$\theta$. Most significantly, we explicitly obtain the following expression for the variance of $\theta$ by comparing
coefficients of $z^2$, which we show (in Proposition \ref{prop10}) are given by the following double integral:

\footnotesize
\beqn
\int_0^{\infty} du_1\int_0^{u_1}\sqrt{\frac{u_1u_2}{\text{sinh}\, u_1 \text{sinh}\, u_2}}\frac{u_1u_2}{u_1+u_2} \parens{\frac{\frac{1}{u_1^2}\parens{1-\frac{u_1}{\text{sinh} u_1}\text{cosh}\,u_1}-\frac{1}{u_2^2}\parens{1-\frac{u_2}{\text{sinh}\, u_2} \text{cosh}\,u_2}}{u_1-u_2}} du_2 .
\eeqn
\normalsize
Although it is not possible to calculate the double integral above in elementary
terms, it has removable singularities and shockingly converges very nicely at all
points where any of $u_1,u_2$, or $u_1-u_2$ vanishes.

\section{A few results needed for obtaining the distribution of $\theta$} \label{sec:formal_statement_results}
In \S \ref{sec21}, we rewrite $\theta$ in an alternate form that will be useful in \S \ref{sec3}. The alternate form involves stochastic integrals rather than integrals of a Wiener process itself. In \S \ref{sec22}, we introduce the function $F$, which is well suited to calculating the moments of $\theta$. Furthermore, we explicitly calculate $F$.  
\subsection{Rewriting $\theta$} \label{sec21}
The author of \cite{magnus} gave the moments of the ratio of a {\em pair of
quadratic forms} in normal variables. Our work
gives a method for the correlation coefficient which is a ratio involving 
{\em three} quadratic forms of normal variables as well as square roots. The three-form problem solved here requires a new method. We begin this task by rewriting $\theta$. 

Recall the definition of empirical correlation written in Equation (\ref{corr}). Noting that $m_i := \int_0^1 W_i(t) dt, i = 1,2$ are the empirical mean values, we rewrite $\theta$ as:

\begin{equation} \label{schwarz}
\theta = \frac{\int_0^1 (W_1(t) - m_1)(W_2(t) - m_2) dt}{\sqrt{\int_0^1(W_1(t) - m_1)^2 dt} \sqrt{\int_0^1 (W_2(t) - m_2)^2 dt}}.
\end{equation}
From Equation ($\ref{schwarz}$) it is easy to see, by Cauchy-Schwarz, that $-1 \le \theta \le 1$. 

\begin{proposition} \label{prop1}
We have the equality
\begin{equation}
\theta = \frac{X_{1,2}}{\sqrt{X_{1,1}X_{2,2}}},
\end{equation}
where
\begin{equation} \label{neweq2}
X_{i,j} = \int_0^1 \int_0^1 \parens{\mathrm{min}(s_1,s_2) -s_1 s_2} dW_i(s_1) dW_j(s_2).
\end{equation}
\end{proposition}

\begin{proof}
It is clear that

\begin{equation*}
\theta = \frac{Y_{1,2}}{\sqrt{Y_{1,1}Y_{2,2}}},
\end{equation*}
where
\begin{equation}
Y_{1,2} = Y_{2,1} = \int_0^1W_1(t)W_2(t) dt - \int_0^1W_1(s) ds\int_0^1 W_2(t) dt
\end{equation}
and
\begin{equation}
Y_{i,i} = \int_0^1 W^2_i(t) dt - \parens{\int_0^1W_i(t) dt}^2, \, i = 1,2.
\end{equation}
We must show that
\beqn
Y_{i,j}=X_{i,j}\,\, \text{for} \,\, i, j \in \{1,2\}.
\eeqn
By the Fundamental Theorem of Calculus, 

\small

\begin{equation} \label{eqcool}
X_{i,j}=\quad \int_0^1 dW_i(s_1) \int_0^1 dW_j(s_2)\int_0^{\text{min}(s_1,s_2)} dt -\int_0^1 dW_i(s_1) \int_0^1 dW_j(s_2) \int_0^{s_1} ds \int_0^{s_2} dt. 
\end{equation}

\normalsize
\noindent To simplify the right hand side of Equation (\ref{eqcool}), we change variables and integrate over the equivalent regions. The right hand side of Equation (\ref{eqcool}) becomes

\small
\begin{equation} \label{newone}
	\quad \,\,\,\,\,\, \int_0^1 dt \int_t^1 dW_i(s_1)\int_t^1 dW_j(s_2) - \int_0^1 ds \int_s^1 dW_i(s_1) \int_0^1dt \int_t^1 dW_j(s_2).
\end{equation}
\normalsize
Rearranging, we obtain

\footnotesize
\beqn
X_{i,j}&=&\int_0^1 \parens{W_i(1)-W_i(t)}\parens{W_j(1) -W_j(t)}dt -\int_0^1 (W_i(1)-W_i(s))ds \int_0^1 (W_j(1) - W_j(t)) dt \nonumber\\
&=& \int_0^1 W_i(t) W_j(t) dt - \int_0^1W_i(s) ds \int_0^1 W_j(t) dt = Y_{i,j}. 
\nonumber
\eeqn

\end{proof}

\normalsize

\subsection{Defining and calculating $F$} \label{sec22}
We define for $|a| \le 1, \beta_i \ge 0, i = 1,2$, the integral

\begin{equation} \label{mgf}
F(\beta_1,\beta_2,a) = \expe{e^{a \beta_1 \beta_2 X_{1,2} - \frac{\beta^2_1}{2} X_{1,1} - \frac{\beta^2_2}{2} X_{2,2}}},
\end{equation}
where the $X_{i,j}$ are as defined in Equation (\ref{neweq2}).
Under the above conditions, the expectation is finite and thus $F$ is well defined. This
is because $X_{1,2} = \theta \sqrt{X_{1,1} X_{2,2} }$ and $|\theta| \le 1$ so
the exponent is at most

\begin{equation}
-\frac{1}{2} \parens{\beta_1 \sqrt{X_{1,1}} - \beta_2 \sqrt{X_{2,2}}}^2 \le 0.
\end{equation}
Thus the expectand is bounded by unity, and so, for this range,
$F(\beta_1,\beta_2,a) \le 1$.\\

\subsubsection{The functions $M$, $K$, $T_{M}$, and $T_{K}$ }
Motivated by the definition of $X_{i,j}$ in Equation (\ref{neweq2}), we define $M$ by
\begin{equation}
M(s_1,s_2) = \text{min}(s_1,s_2) - s_1 s_2, \,\, s_1,s_2 \in [0,1],
\end{equation}
which is the covariance of pinned Wiener process on $[0,1]$.

For $i_1,i_2 \in \{1,2\}$ and $s_1,s_2 \in [0,1]$, we define the kernel function $K_{i_1,i_2}$ by

\begin{equation} \label{K1}
K_{i,i}\parens{s_1,s_2}= -\beta^2_i M\parens{s_1,s_2}
\end{equation}
and
\begin{equation} \label{K2}
K_{1,2}(s_1,s_2) = K_{2,1}(s_1,s_2) =a\beta_1 \beta_2 M(s_1,s_2). \nonumber \\
\end{equation}
Let
\beqn
K=\twobytwomat{K_{1,1}}{K_{1,2}}{K_{2,1}}{K_{2,2}}.
\eeqn
$M$ gives rise to a linear transformation $T_M: L^2[0,1]\rightarrow L^2[0,1]$ given by
\begin{equation}
T_M(g)(s_2)=\int_0^1 M(s_1,s_2)g(s_1)ds_1
\end{equation}
where $g \in L^2[0,1]$. If $\vec{f}(s)$ and $\vec{g}(s)$ are elements of $\parens{L^2[0,1]}^2$, then their inner product is defined as
\beqn
\int_0^1 \vec{f}(s)\cdot \vec{g}(s)ds.
\eeqn          It is straightforward to check that $T_M$ is a self-adjoint, positive-definite, and continuous linear transformation.

Let $M_{2\times 2}$ refer to 2 $\times$ 2 matrices with entries in $L^2([0,1]\times [0,1])$. We define $\parens{L^2[0,1]}^2=L^2[0,1] \times L^2[0,1]$ and note that an element of $\parens{L^2[0,1]}^2$ is an ordered pair of maps from $[0,1]$ to $\mathbb{R}$. If $G$ is any element of $M_{2 \times 2}(L^2([0,1]\times [0,1]))$, $G$ gives rise to a linear transformation $T_G: \parens{L^2[0,1]}^2\rightarrow \parens{L^2[0,1]}^2$ given by
\begin{equation} \label{TG}
T_G(\vec{g})(s_2)=\int_0^1 G(s_1,s_2)\vec{g}(s_1)ds_1
\end{equation}
for $\vec{g} \in \parens{L^2[0,1]}^2$. 
In particular, $K$ gives rise to $T_K$ via the above definition. It is straightforward to check that $T_K$ is a self-adjoint, positive-definite, and continuous linear transformation.

\subsubsection{The calculation}

The goal of this section is to prove Theorem \ref{nicetheorem} below.

\begin{theorem} \label{nicetheorem}
\begin{equation}
F\parens{\beta_1,\beta_2,a}= \frac{1}{\sqrt{\frac{\sinh{c^+}}{c^+} \frac{\sinh{c^-}}{c^-}}}, \label{eqn:F}
\end{equation}
where 
\begin{equation} \label{eq33}
\quad \,\,\,\,\,\,c^\pm = c^\pm(\beta_1,\beta_2,a) = \sqrt{\frac{(\beta^2_1 + \beta^2_2) \pm \sqrt{(\beta^2_1-\beta^2_2)^2 +4a^2 \beta^2_1 \beta^2_2}}{2}}.
\end{equation}
\end{theorem}

The proof of Theorem \ref{nicetheorem} has two parts. First, we show (in Proposition \ref{prop4}) that $F$ can be expressed in terms of the Fredholm determinant $\text{det}(I-T_K)$, where $T_K$ was defined in Equation (\ref{TG}). Second, we compute this determinant (in Proposition \ref{prop6}) by calculating the eigenvalues of $T_K$. To begin, we write $F$ in terms of $K$.
\begin{proposition} 
We have
\small
\begin{equation} \label{eqn3}
\quad \,\,\,\,\,\, F(\beta_1,\beta_2,a) = \expe{e^{\frac{1}{2}\sum_{i_1 =1}^2 \sum_{i_2 = 1}^2 \int_0^1 \int_0^1 K_{i_1,i_2}\parens{s_1,s_2} dW_{i_1}(s_1) dW_{i_2}(s_2)}}.
\end{equation}
\end{proposition}
\normalsize
\begin{proof}
This follows from plugging in the definition of $K_{i_1,i_2}$ from Equations (\ref{K1}) and (\ref{K2}) and comparing it to the definition of $F$ given in Equation (\ref{mgf}).
\end{proof}
Since $T_K$ is self-adjoint and positive-definite, there exists a (countable) orthonormal basis of $(L^2[0,1])^2$ consisting of eigenvectors $\vec{\phi}_n\,  (n \in \mathbb{N})$ for $T_K$. As a matter of notation, we write
\beqn
\vec{\phi}_n= \twovec{\phi_n(1,s_1)}{\phi_n(2,s_1)}.
\eeqn
We let $\alpha_n$ be the corresponding eigenvalues. 
\begin{proposition}\label{Kprop}
For $i_1,i_2 \in \{1,2\}, s_1,s_2 \in [0,1],$ we have
\begin{equation} \label{newone}
K_{i_1,i_2}(s_1,s_2) = \sum_{n = 1}^\infty \alpha_n \phi_n(i_1,s_1)\phi_n(i_2,s_2).
\end{equation}
\end{proposition}
\begin{proof}
Let $G_n(s_1,s_2)= \vec{\phi}_n \vec{\phi}_n^{\,T} \in M_{2 \times 2}(L^2([0,1] \times [0,1]))$. The proposition is equivalent to the equation $K=\sum_{n=1}^\infty \alpha_nG_n$. Thus, we need only show that
\begin{equation} \label{coolneweq}
T_K=\sum_{n=1}^\infty \alpha_n T_{G_n}.
\end{equation}
Since the vectors
$\vec{\phi}_n$ form a basis for $\parens{L^2[0,1]}^2$, it suffices to check both sides of Equation (\ref{coolneweq}) on an arbitrary basis element $\vec{\phi}_m$.

By definition, the left hand side of Equation (\ref{coolneweq}) applied to $\vec{\phi}_m$ yields $\alpha_m \vec{\phi}_m$ (note that $\vec{\phi_m}$ is chosen among the eigenvectors for $T_K$, i.e.,  it is one of the $\vec{\phi_n}$). So it suffices to show that
\beqn
T_{G_n}(\vec{\phi}_m)= \delta_{m,n}\vec{\phi}_m.
\eeqn
We calculate
\beqn
T_{G_n}(\vec{\phi}_m)&=&\int_0^1  \twovec{\phi_n(1,s_2)}{\phi_n(2,s_2)} \twovec{\phi_n(1,s_1)}{\phi_n(2,s_1)}^T \twovec{\phi_m(1,s_1)}{\phi_m(2,s_1)}ds_1\\
&=& \twovec{\phi_n(1,s_2)}{\phi_n(2,s_2)} \int_0^1 \vec{\phi}_n \cdot \vec{\phi}_m \, ds_1.
\eeqn
Since $ \vec{\phi}_n$ and $ \vec{\phi}_m$ are orthonormal, 
\beqn
T_{G_n}(\vec{\phi}_m)= \twovec{\phi_n(1,s_2)}{\phi_n(2,s_2)} \delta_{m,n}= \delta_{m,n}\vec{\phi}_m.
\eeqn
The proof is now complete.
\end{proof}

\begin{proposition} \label{prop4}
$F(\beta_1, \beta_2, a)= \parens{\mathrm{det}(I-T_K)}^{-1/2}:= \prod_n \parens{1-\alpha_n}^{-1/2}$.
\end{proposition}
\begin{proof}
By substituting Equation (\ref{newone}) into Equation (\ref{eqn3}), we obtain
\begin{equation}\label{eqn:F}
F\parens{\beta_1,\beta_2,a}= \prod_n \expe{e^{\frac{1}{2}\parens{\sum_{i = 1}^2 \alpha_n \int_0^1 \phi_n(i,s_i) dW(s_{i})}^2}}.
\end{equation}
Letting $\xi_i$ be real numbers, the right hand side of Equation (\ref{eqn:F}) is equal to $\prod_n \parens{1-\alpha_n}^{-1/2}$. 
\end{proof}
Now, we proceed to calculate the eigenvalues of $T_K$. We do this by ``guessing'' the eigenvectors. Because the entries in $K$ are scalar multiples of each other, we guess that the eigenvectors are of separable
form:
\begin{equation}
\vec{\phi}(s) = \twovec{\xi_{1} \phi(s)}{\xi_{2} \phi(s)},\,\,s \in [0,1],
\end{equation}
where $\phi(s) \in L^2[0,1]$ is an eigenvector of $T_M$ and $\xi_i$ are real numbers. It can be easily shown that these ``guessed'' eigenvectors span the entire $\parens{L^2[0,1]}^2$.

It is straightforward to verify that the functions $\psi_n(t):=\sqrt{2} \sin{\pi n t}$ for $n \geq 1$ and $t \in [0,1]$ form an orthonormal basis of $L^2[0,1]$ consisting of eigenvectors for $T_M$.  The eigenvalue $\lambda_n$ corresponding to $\psi_n$ can be calculated to be $\frac{1}{\pi^2 n^2}$, i.e.
\begin{equation*} \label{clarification}
T_M\psi_n(t)=\frac{1}{\pi^2 n^2}\psi_n(t).
\end{equation*}
We find that $\twovec{\xi_{1} \phi(s)}{\xi_{2} \phi(s)}$ is an eigenvector of $T_K$ with eigenvalue $\alpha$ if and only if

\begin{equation} \label{sep1}
-\beta^2_1 \xi_1 + a \beta_1 \beta_2 \xi_2 = \frac{\alpha}{\lambda_n} \xi_1
\end{equation}
and

\begin{equation} \label{sep2}
a \beta_1 \beta_2 \xi_1  -\beta^2_2 \xi_2 = \frac{\alpha}{\lambda_n} \xi_2.
\end{equation} 

\begin{proposition} \label{prop5}
For each eigenvalue $\lambda_n$ of $T_M$, there are two corresponding eigenvalues $\gamma_n^{\pm}$ of $T_K$, where 
\beqn
\gamma_n^{\pm}= \lambda_n \frac{-(\beta^2_1 + \beta^2_2) \pm \sqrt{(\beta^2_1-\beta^2_2)^2 +4a^2 \beta^2_1 \beta^2_2}}{2}.
\eeqn
\end{proposition}
\begin{proof}
To find the eigenvalues, we solve Equations (\ref{sep1}) and (\ref{sep2}) for $\alpha_n$. We can view the system of Equations (\ref{sep1}) and (\ref{sep2}) as the system $A\vec{v}=\lambda \vec{v}$, where 
\beqn
A=\twobytwomat{-\beta_1^2}{a \beta_1 \beta_2}{a \beta_1 \beta_2}{-\beta_2^2}
\eeqn
and 
\beqn
\vec{v}=\twovec{\xi_1}{\xi_2}
\eeqn
and
\beqn
\lambda=\frac{\alpha}{\lambda_n}.
\eeqn
The solutions for $\alpha$ are $\lambda_n$ multiplied by the eigenvalues of $A$. Calculating the eigenvalues of $A$ finishes the proposition.
\end{proof}

\begin{proposition} \label{prop6}
If $\alpha_n, \, n\geq 1$ is the list of eigenvalues for $T_K$, then
\beqn
\prod_{n=1}^\infty (1-\alpha_n)=\frac{\sinh{c^+}}{c^+} \frac{\sinh{c^-}}{c^-},
\eeqn
where $c^\pm$ is as defined in Equation (\ref{eq33}).
\end{proposition}
\begin{proof}
First note that the set \{$\alpha_n, \, n \geq 1$\} is the same as the set \{$\gamma^{\pm}_n,\, n \geq 1$\}. Thus, 
\bneqn \label{boaseq}
\prod_{n=1}^\infty (1-\alpha_n)&=&\prod_{n=1}^\infty (1-\gamma_n^+) (1-\gamma_n^-) \nonumber \\
&=& \prod_{n=1}^\infty \parens{1 - \frac{(z^+)^2}{\pi^2 n^2}} \prod_{n=1}^\infty \parens{1- \frac{(z^-)^2}{\pi^2 n^2}}= \frac{\sin{z^+}}{z^+} \frac{\sin{z^-}}{z^-},
\eneqn
where
\begin{equation} 
z^\pm = \sqrt{\frac{-(\beta^2_1 + \beta^2_2) \pm \sqrt{(\beta^2_1-\beta^2_2)^2 +4a^2 \beta^2_1 \beta^2_2}}{2}}.
\end{equation}

\end{proof}
The final equality in Equation (\ref{boaseq}) follows from the following product formula (see \cite{boas})

\begin{equation}
\frac{\sin{z}}{z} = \prod_{n=1}^\infty \parens{1 - \frac{z^2}{\pi^2 n^2}}.
\end{equation}
Note that for any complex number $z$,
\beqn
\frac{\sin{z}}{z} = \frac{\sinh{(-iz)}}{-iz}.
\eeqn
Since $z^\pm$ are purely imaginary, we observe $z^{\pm} = ic^\mp$ and the $c^\mp$ are nonnegative real. In particular, since $c^\mp=-iz^{\pm}$, we have that 
\beqn
\frac{\sin{z^{+}}}{z^{+}}\frac{\sin{z^{-}}}{z^{-}} = \frac{\sinh{(-iz^{+})}}{-iz^{+}}\frac{\sinh{(-iz^{-})}}{-iz^{-}}=\frac{\sinh{c^{+}}}{c^{+}}\frac{\sinh{c^{-}}}{c^{-}}.
\eeqn
This completes the proof.\\

\noindent \textbf{Proof of Theorem \ref{nicetheorem}}. The proof immediately follows from combining Proposition \ref{prop4} and Proposition \ref{prop6}. \qed

\begin{remark}
Both of $c^\pm$ are nonnegative if $|a| \le 1$, but if 

\begin{equation}
|a| \ge \sqrt{1+\frac{2(\beta^2_1+\beta^2_2)\pi^2 + \pi^4}{4\beta^2_1 \beta^2_2}}
\end{equation}
then $F(\beta_1,\beta_2,a) = \infty$ since the term involving $c^-$ vanishes
as $c^-$ becomes imaginary and the term in the denominator becomes
$\sin{-\pi} = 0$. However, we will only need $F(\beta_1,\beta_2,a)$ for
$|a| < 1$ to obtain the distribution of $\theta$.

\end{remark}

\section{Obtaining the integral equation for the moment generating function of $\theta$} \label{sec3}
In this section, we derive a formula for the moments of $\theta$. Our strategy is to use integral equations to derive the distribution of $\theta$. The methods we employ draw inspiration from some related ideas and approaches in earlier works of \cite{shepp1,shepp2,shepp3, shepp4}.  

\subsection{The Form for the Moments of $\theta$} \label{sec31}

In this chapter, we will generally denote derivatives with subscripts. However, we will sometimes denote derivatives of functions of three arguments with respect to the third argument by primes rather than subscripts because it is more natural for Theorem \ref{thm2} below. For example, 

\small
\beqn
F_i(\beta_1,\beta_2,z) = \frac{\partial}{\partial \beta_i}F(\beta_1,\beta_2,z),
i = 1,2; \,\,\, F_3(\beta_1,\beta_2,z) =F^\prime(\beta_1,\beta_2,z) = \frac{\partial}{\partial z}F(\beta_1,\beta_2,z).
\eeqn
\normalsize
The goal of \S \ref{sec31} is to prove the following:
\begin{theorem} \label{thm2}
The moments of $\theta$, where $\theta$ is defined in Equation (\ref{corr}),  satisfy
\small
\beqn
\sum_{n=1}^\infty \frac{z^{2n}}{2n}  \expe{\theta^{2n}} \frac{(n!)^2 2^{2n}}{(2n)!} = \int_0^\infty \frac{d\beta_1}{\beta_1} \int_0^\infty \frac{d\beta_2}{\beta_2} z F^\prime(\beta_1,\beta_2,z) \label{eqn:simpler}.
\eeqn
\end{theorem}

We first introduce the function
$G = G(\gamma_1,\gamma_2,a)$, given by

\begin{equation} \label{construction}
G=\int_0^\infty \frac{d \beta_1}{\beta_1} \int_0^\infty \frac{d \beta_2}{\beta_2}\bar{V},
\end{equation}
where $\bar{V}$ is defined as:

\footnotesize
\beqn
F(\beta_1,\beta_2,a) - F\parens{\sqrt{\gamma_1}\beta_1,\beta_2,\frac{a}{\sqrt{\gamma_1}}} - F\parens{\beta_1,\sqrt{\gamma_2}\beta_2,\frac{a}{\sqrt{\gamma_2}}} + F\parens{\sqrt{\gamma_1}\beta_1,\sqrt{\gamma_2}\beta_2,\frac{a}{\sqrt{\gamma_1 \gamma_2}}}.
\eeqn
\normalsize

\begin{remark} \label{conditionsremark}
Note that if $\gamma_i \ge 1, \beta_i \ge 0, i = 1,2$, and $|a| < 1$, the quantities 
\beqn
F\parens{\sqrt{\gamma_1}\beta_1,\beta_2,\frac{a}{\sqrt{\gamma_1}}}, F\parens{\beta_1,\sqrt{\gamma_2}\beta_2,\frac{a}{\sqrt{\gamma_2}}}, F\parens{\sqrt{\gamma_1}\beta_1,\sqrt{\gamma_2}\beta_2,\frac{a}{\sqrt{\gamma_1 \gamma_2}}}
\eeqn
are finite and well defined. 
One can (with some algebra) show that the right hand side of Equation (\ref{construction}) converges and thus $G$ is also well defined.
\end{remark}

\begin{proposition} \label{propderivatives}
Under the conditions of Remark $\ref{conditionsremark}$, the expression
\beqn
\frac{\partial^2 G(\gamma_1,\gamma_2,a)}{\partial \gamma_1 \partial \gamma_2}
\eeqn
is simultaneously equal to

\small
\begin{equation} \label{coolone}
\int_0^\infty d\beta_1 \int_0^\infty d\beta_2 \frac{\beta_1 \beta_2}{4} \expe{e^{a \beta_1 \beta_2 \theta} e^{-\gamma_1 \frac{\beta^2_1}{2}} e^{- \gamma_2 \frac{\beta^2_2}{2} }}
\end{equation}
\normalsize
and
\small
\begin{equation} \label{neweq}
\int_0^\infty \frac{d \beta_1}{\beta_1} \int_0^\infty \frac{d \beta_2}{\beta_2} \frac{\partial^2}{\partial \gamma_1 \partial \gamma_2} F\parens{\sqrt{\gamma_1} \beta_1,\sqrt{\gamma_2} \beta_2,\frac{a}{\sqrt{\gamma_1 \gamma_2}}}.
\end{equation}
\end{proposition}

\begin{proof}
Plugging in the definition of $F$ from Equation (\ref{mgf}), Equation ($\ref{construction}$) can be rewritten as

\footnotesize
\begin{equation*}
G(\gamma_1,\gamma_2,a) = \int_0^\infty \frac{d \beta_1}{\beta_1} \int_0^\infty \frac{d \beta_2}{\beta_2} \expe{e^{a\beta_1\beta_2 X_{1,2}} \parens{e^{-\frac{\beta^2_1 X_{1,1}}{2}} - e^{-\frac{\gamma_1 \beta^2_1 X_{1,1}}{2}}} \parens{e^{-\frac{\beta^2_2 X_{2,2}}{2}} - e^{-\frac{\gamma_2 \beta^2_2 X_{2,2}}{2}}}}.
\end{equation*}
\normalsize
The second integrand of the above equation is everywhere positive and if we take the
expectation outside the integrals, then for each fixed
$\omega \in \Omega$ (the probability space where $W_i(s,\omega)$ are defined), we can replace $\beta_i$,
by $\frac{\beta_i}{\sqrt{X_{i,i}(\omega)}}, i = 1,2$. We then rewrite $G(\gamma_1,\gamma_2,a)$ as follows:

\small
\beqn
G(\gamma_1,\gamma_2,a) = \expe{\int_0^\infty \frac{d \beta_1}{\beta_1} \int_0^\infty \frac{d \beta_2}{\beta_2} e^{a \beta_1 \beta_2 \theta(\omega)} \parens{e^{-\frac{\beta^2_1}{2}} - e^{-\frac{\gamma_1 \beta^2_1}{2}}} \parens{e^{-\frac{\beta^2_2}{2}} - e^{-\frac{\gamma_2 \beta^2_2}{2}}}}.
\eeqn
\normalsize
Putting the expectation back inside the integral, we obtain

\footnotesize
\begin{equation}\label{eqn:key} 
G(\gamma_1,\gamma_2,a) = \int_0^\infty \frac{d \beta_1}{\beta_1} \int_0^\infty \frac{d \beta_2}{\beta_2} \expe{e^{a \beta_1 \beta_2 \theta}} \parens{e^{-\frac{\beta^2_1}{2}} - e^{-\frac{\gamma_1 \beta^2_1}{2}}} \parens{e^{-\frac{\beta^2_2}{2}} - e^{-\frac{\gamma_2 \beta^2_2}{2}}}.
\end{equation}
\normalsize


\normalsize


Recall the corollary of Fubini's theorem that states that integration and
differentiation with respect to a parameter can be interchanged if the integral of the
differentiated integrand converges absolutely. One can show this is true in Equation (\ref{eqn:key}) for $\abss{a}<1$.
Thus, for $\abss{a}<1$, we obtain

\small
\beqn 
\frac{\partial^2 G(\gamma_1,\gamma_2,a)}{\partial \gamma_1 \partial \gamma_2}=\int_0^\infty d\beta_1 \int_0^\infty d\beta_2 \frac{\beta_1 \beta_2}{4} \expe{e^{a \beta_1 \beta_2 \theta}} e^{-\gamma_1 \frac{\beta^2_1}{2}} e^{- \gamma_2 \frac{\beta^2_2}{2} }.
\eeqn
\normalsize
On the other hand, using the definition of $G$ from Equation (\ref{construction}) directly, we obtain

\small
\beqn 
\,\,\quad  \frac{\partial^2 G(\gamma_1,\gamma_2,a)}{\partial \gamma_1 \partial \gamma_2} &=& \int_0^\infty \frac{d \beta_1}{\beta_1} \int_0^\infty \frac{d \beta_2}{\beta_2} \frac{\partial^2}{\partial \gamma_1 \partial \gamma_2} F\parens{\sqrt{\gamma_1} \beta_1,\sqrt{\gamma_2} \beta_2,\frac{a}{\sqrt{\gamma_1 \gamma_2}}}.
\eeqn
\normalsize
This completes the proof.
\end{proof}

\normalsize
Our next proposition gives our first expression for the moments of $\theta$. We now fix some notation. Let $F_i$ denote the derivative with respect to the $i$th argument, as follows:
\beqn
F_i(\beta_1, \beta_2, z)&=& \frac{\partial}{\partial \beta_i}F(\beta_1, \beta_2, z), \,i=1,2;\\
F_3(\beta_1, \beta_2, z)&=& F'(\beta_1, \beta_2, z)= \frac{\partial}{\partial z}F(\beta_1, \beta_2, z).
\eeqn
\begin{proposition} \label{prop8}
We have

\footnotesize
\bneqn \label{brackseq}
&&\sum_{n=0}^\infty z^{2n} \expe{\theta^{2n}} \frac{(n!)^2 2^{2n}}{(2n)!} \nonumber \\
&=&\int_0^\infty \frac{d\beta_1}{\beta_1} \int_0^\infty \frac{d\beta_2}{\beta_2} \parens{z^2 F^{\prime \prime}(\beta_1, \beta_2, z) + z F^\prime(\beta_1, \beta_2, z) - z \beta_1 F^\prime_1(\beta_1, \beta_2, z) - z \beta_2 F^\prime_2(\beta_1, \beta_2, z) + \beta_1 \beta_2 F_{12}(\beta_1, \beta_2, z)} \nonumber
\eneqn
\end{proposition}
\normalsize

\begin{proof}
Our strategy is to simplify expressions (\ref{coolone}) and (\ref{neweq}) and to set them equal to each other using Proposition  \ref{propderivatives}.  We begin now by simplifying the expression in (\ref{coolone}). We do this by expanding $\expe{e^{a\beta_1 \beta_2 \theta}}$ as a power series
and substitute $\frac{\beta^2_i}{2}$ with $u_i$,
$i = 1,2$, and do the integrals on $u_i$, term-by-term, to obtain an
expression involving the moments, $\mu_{2n} = \expe{\theta^{2n}}$, of $\theta$ 
(note that the odd moments, $\mu_{2n+1}, n \ge 0$, vanish, by symmetry). Algebraic manipulation yields a simpler form of the expression in (\ref{coolone}), displayed below:

\small
\bneqn \label{nice}
&& \sum_{n=0}^\infty \frac{a^{2n}}{(2n)!} \mu_{2n} \frac{1}{4} \int_0^\infty d u_1 \int_0^\infty  (2u_1)^n e^{-\gamma_1 u_1} (2u_2)^n e^{-\gamma_1 u_2} d u_2\nonumber\\
&=& \sum_{n=0}^\infty \parens{\frac{a}{\sqrt{\gamma_1 \gamma_2}}}^{2n} \mu_{2n} \frac{(n!)^2 2^{2n}}{(2n)!} \frac{1}{4\gamma_1 \gamma_2}.
\eneqn
\normalsize

Thus we see that the simplified form of expression (\ref{coolone}) is $\frac{1}{4\gamma_1 \gamma_2}$ times a
function of $z = \frac{a}{\sqrt{\gamma_1 \gamma_2}}$, and of course
the expression (\ref{neweq}) must also be of this form.
\normalsize
Applying the multivariate chain rule to the expression $\frac{\partial^2}{\partial \gamma_1 \partial \gamma_2} F\parens{\sqrt{\gamma_1} \beta_1,\sqrt{\gamma_2}\beta_2,\frac{a}{\sqrt{\gamma_1 \gamma_2}}}$ from expression (\ref{neweq}) yields

\small
\beqn 
\frac{\partial}{\partial \gamma_2}\parens{\frac{1}{2} \frac{\beta_1}{\sqrt{\gamma_1}} F_1\parens{\sqrt{\gamma_1}\beta_1,\sqrt{\gamma_2}\beta_2,\frac{a}{\sqrt{\gamma_1 \gamma_2}}} + \frac{-a}{2 \gamma_1^\frac{3}{2} \gamma_2^\frac{1}{2}} F_3\parens{\sqrt{\gamma_1}\beta_1,\sqrt{\gamma_2}\beta_2,\frac{a}{\sqrt{\gamma_1 \gamma_2}}}}
\eeqn
\normalsize
which in turn equals

\scriptsize
\begin{equation} \label{none}
\quad \quad \frac{\beta_1 \beta_2}{4 \sqrt{\gamma_1 \gamma_2}} F_{12} + \frac{-a \beta_1}{4\sqrt{(\gamma_1)^2 (\gamma_2)^3}} F_{13} + \frac{-a \beta_2}{4\sqrt{(\gamma_1)^3 (\gamma_2)^2}} F_{23} + \frac{a}{4 \sqrt{(\gamma_1)^3 (\gamma_2)^3}} F_3 + \frac{a^2}{4 \gamma^2_1 \gamma^2_2} F_{33}.
\end{equation}
\normalsize
Note that in expression (\ref{none}), for the sake of brevity, we dropped the arguments

\begin{equation}
\parens{\sqrt{\gamma_1}\beta_1,\sqrt{\gamma_2}\beta_2,\frac{a}{\sqrt{\gamma_1 \gamma_2}}}
\end{equation}
of $F$ and its derivatives. Substituting the right hand side of expression (\ref{none}) into expression (\ref{neweq}) and replacing $\beta_i$ by $\frac{\beta_i}{\sqrt{\gamma_i}}, i = 1,2$, and setting $z = \frac{a}{\sqrt{\gamma_1 \gamma_2}}$,  we obtain

\small
\begin{equation} \label{eq40}
\quad \frac{1}{4\gamma_1 \gamma_2} \int_0^\infty \frac{d\beta_1}{\beta_1} \int_0^\infty \frac{d\beta_2}{\beta_2} \parens{z^2 F^{\prime \prime} + z F^\prime - z \beta_1 F^\prime_1 - z \beta_2 F^\prime_2 + \beta_1 \beta_2 F_{12}}.
\end{equation}
\normalsize

The proposition now follows by equating expression (\ref{eq40}) with the expression in (\ref{nice}). 
\end{proof}

The following proposition completes the proof of Theorem \ref{thm2}.
\normalsize
\begin{proposition} \label{prop9}
Equation (\ref{brackseq}) can be simplified to 
\begin{equation} \label{eqn:simpler}
\sum_{n=1}^\infty \frac{z^{2n}}{2n}  \expe{\theta^{2n}} \frac{(n!)^2 2^{2n}}{(2n)!} = \int_0^\infty \frac{d\beta_1}{\beta_1} \int_0^\infty \frac{d\beta_2}{\beta_2} z F^\prime(\beta_1,\beta_2,z).
\end{equation}
\end{proposition}
\normalsize

\begin{proof}

We first consider the third, fourth, and fifth terms in the integrand of the right hand side of Equation (\ref{brackseq}). The third and fourth terms give zero upon integration. For
example, the fourth term is again an exact differential in $\beta_2$, so

\begin{equation}
\int_0^\infty d\beta_2-\int_0^\infty \frac{d\beta_1}{\beta_1} zF^\prime_2(\beta_1,\beta_2,z) = z \int_0^\infty \frac{d\beta_1}{\beta_1} F^\prime(\beta_1,0,z) = 0
\end{equation}
because $F(\beta_1,0,z) = \expe{e^{-\frac{\beta^2_2}{2}X_{1,1}}}$, which is constant
in $z$ and so has a derivative of zero with respect to $z$. The same argument holds for the third term.
The fifth term gives unity upon integration because the $\beta_1 \beta_2$ terms in
the numerator and denominator cancel and the integrand becomes an exact differential, giving

\small
\beqn
\int_0^\infty d\beta_1 \int_0^\infty d\beta_2 F_{12}(\beta_1,\beta_2,z) = - \int_0^\infty d \beta_1 F_1(\beta_1,0,z) = F(0,0,z) = Ee^{0} = 1.
\eeqn
\normalsize
The sum of the first and second terms can be written as:
\begin{equation} z^2F^{\prime \prime}(\beta_1,\beta_2,z) + z F^\prime(\beta_1,\beta_2,z) = z\frac{d}{dz} \parens{zF^\prime(\beta_1,\beta_2,z)}.
\end{equation}
We may subtract unity from both sides of Equation (\ref{brackseq}), divide by $z$, and
integrate to obtain

\begin{equation} \label{eqn:simpler}
\sum_{n=1}^\infty \frac{z^{2n}}{2n}  \expe{\theta^{2n}} \frac{(n!)^2 2^{2n}}{(2n)!} = \int_0^\infty \frac{d\beta_1}{\beta_1} \int_0^\infty \frac{d\beta_2}{\beta_2} z F^\prime(\beta_1,\beta_2,z).
\end{equation}
This finishes the proof.
\end{proof}

\subsection{Isolating the Moments}
The goal of this section is to use Proposition \ref{prop9} to extract the even moments of $\theta$. We begin by defining the function $S(u)$ as:

\beqn
S(u)=\sqrt{\frac{u}{\text{sinh}\, u}}.
\eeqn
Recalling the definition of $c^{\pm}$ in Equation (\ref{eq33}), Theorem \ref{nicetheorem} says

\beqn
F(\beta_1,\beta_2,z) = S(c^+(\beta_1,\beta_2,z) ) S(c^-(\beta_1,\beta_2,z) ).
\eeqn
Differentiating with respect to $z$ and multiplying $z$ yields

\beqn
z \frac{d}{dz} F(\beta_1,\beta_2,z) = z S^\prime(c^+) (c^+)^\prime S(c^-) + z S(c^+) S^\prime(c^-) (c^-)^\prime.
\eeqn
(note that, as before, the prime notation denotes the derivative with respect to $z$). For example, the derivative of $c^{+}$ (with respect to $z$) is, after some algebraic manipulation, simplified as:
\beqn
\frac{1}{\sqrt{c^+}}\frac{z b_1^2 b_2^2}{\sqrt{(b_1^2-b_2^2)^2+ 4z^2 b_1^2 b_2^2}}.
\eeqn

Define the function $T(c)$ as follows:

\begin{equation}
T(c) = \frac{1}{c} \frac{S^\prime(c)}{S(c)} = \frac{1}{2c^2}\parens{1 - \frac{c}{\sinh{c}} \cosh{c}}.
\end{equation}
Then Equation (\ref{eqn:simpler}) becomes
\footnotesize
\begin{equation} \label{symmetriceq}
\sum_{n=1}^\infty \frac{z^{2n}}{2n} \expe{\theta^{2n}} \frac{(n!)^2 2^{2n}}{(2n)!} = z^2 \int_0^\infty \int_0^\infty S(c^+)S(c^-) \parens{T(c^+) - T(c^-)} \frac{\beta_1 \beta_2 d\beta_1 d\beta_2}{\sqrt{(\beta^2_1-\beta^2_2)^2 + 4\beta^2_1 \beta^2_2 z^2}}.
\end{equation}
\normalsize
We simplify Equation (\ref{symmetriceq}) by noting that we can break up the integrand by finding two regions such that the integral is the same over each of them. Noting that the integral is symmetric in $(\beta_1,\beta_2)$, we simplify Equation (\ref{symmetriceq}) as
\small
\begin{equation} \label{massive}
\sum_{n=1}^\infty \frac{z^{2(n-1)}}{2n}  \expe{\theta^{2n}} \frac{(n!)^2 2^{2n}}{(2n)!} = 2 \int_0^\infty d\beta_1 \int_0^{\beta_1} \beta_1 \beta_2 S(c^+)S(c^-)\frac{T(c^+) - T(c^-)}{(c^+)^2 - (c^-)^2}d\beta_2 . 
\end{equation}
\normalsize
A change of variables will suffice to simplify Equation (\ref{massive}). We proceed by changing variables from $(\beta_1, \beta_2)$ to $(u_1, u_2)$, where, for fixed real $z \in (0,1)$, we set:

\begin{equation}
u_1(\beta_1, \beta_2)=c^{+}(\beta_1,\beta_2,z), \quad u_2(\beta_1, \beta_2)=c^{-}(\beta_1,\beta_2,z).
\end{equation}
The Jacobian of
this transformation can be calculated as:
$J=\frac{\beta_1^2-\beta_2^2}{u_1^2-u_2^2}\sqrt{1-z^2}$, and for
further reference we note that $\beta_1^2+\beta_2^2=u_1^2+u_2^2$ and
$\beta_1\beta_2\sqrt{1-z^2}=u_1u_2$, as well as $d\beta_1d\beta_2=\frac{1}{J}du_1du_2$. After bringing $\sqrt{1-z^2}$ to the left side and canceling a factor of $u_1^2-u_2^2$, Equation (\ref{massive}) becomes

\small
\begin{equation} \label{massive2}
\sqrt{1-z^2} \sum_{n=1}^\infty \frac{z^{2(n-1)}}{2n} \expe{\theta^{2n}}\frac{(n!)^22^{2n}}{(2n)!}= 2 \iint_U \frac{u_1u_2S(u_1)S(u_2)\parens{T(u_1)-T(u_2)}}{\sqrt{\parens{u_1^2-u_2^2}^2-z^2\parens{u_1^2+u_2^2}^2}}du_1du_2.
\end{equation}
\normalsize
Here, $U$ is the region of $(u_1,u_2)$ in which the quantity in the square root of Equation (\ref{massive2}) is positive, i.e.
\begin{equation} \label{eq318}
U=\left\{(u_1,u_2): 0 \leq u_2 \leq u_1\sqrt{\frac{1-z}{1+z}}\right\}.
\end{equation}
To determine $U$, first note that the domain of integration in Equation (\ref{massive}) is bounded by $\beta_1 = \beta_2$ and $\beta_2 = 0$ in the first quadrant.  We make the following change of variables: $u_1=u\sqrt{1+v}$ and $u_2=u\sqrt{1-v}$. Noting that $u_1 \geq u_2$, the two boundaries become $u_2 = u_1\sqrt{\frac{1-z}{1+z}}$ and $u_2 = 0$, which is the description of $U$ in Equation (\ref{eq318}).  The fact that $\beta_2 = 0$ is equivalent to $u_2 = 0$ follows from $\beta_1\beta_2 = u_1u_2\sqrt{1-z^2}$ along with $z \in (0,1)$, $0 \leq u_2 \leq u_1$, and $0 \leq \beta_2 \leq \beta_1$.  The fact that $\beta_1 = \beta_2$ is equivalent to $u_2 = u_1\sqrt{\frac{1-z}{1+z}}$ can be derived directly from the formulas for $c^+$ and $c^-$ in Equation (\ref{eq33}), with $z$ in place of $a$ in Equation (\ref{eq33}).\\
\indent For $0 \leq z <1$, Equation (\ref{massive2}) becomes

\footnotesize
\begin{equation} \label{niceone}
\sqrt{1-z^2}\sum_{n=1}^\infty \frac{z^{2(n-1)}}{2n}  \expe{\theta^{2n}} \frac{(n!)^22^{2n}}{(2n)!}=
2\int_z^1dv \int_0^\infty \frac{uS(u\sqrt{1+v})S(u\sqrt{1-v})\parens{(T(u\sqrt{1+v})-T(u\sqrt{1-v})}}{\sqrt{v^2-z^2}} du .
\end{equation}
\normalsize
We now seek to simplify Equation (\ref{niceone}). It is easy to check that the right hand side of Equation (\ref{niceone}) converges. Let us define $g(v)$ as:

\footnotesize
\begin{equation} \label{alsonice}
g(v)=\int_0^\infty uS(u\sqrt{1+v})S(u\sqrt{1-v})T(u\sqrt{1+v})du-\int_0^\infty uS(u\sqrt{1+v})S(u\sqrt{1-v})T(u\sqrt{1-v})du.
\end{equation}
\normalsize
By the definition of $S$ and $T$, and with some algebra, it can be shown that $g(v)$ is the difference of convergent integrals. 
\normalsize
We now make the substitution $u\rightarrow \frac{u}{\sqrt{1+v}}$ in the first integral in Equation (\ref{alsonice}) and the substitution $u\rightarrow \frac{u}{\sqrt{1-v}}$ in the second integral in Equation (\ref{alsonice}), we can simplify the right hand side of Equation (\ref{niceone}) as:

\small
\begin{equation} \label{niceone2}
2\int_z^1\frac{dv}{\sqrt{v^2-z^2}}\int_0^\infty S(u)T(u)\bracks{\frac{u}{1+v}S\parens{u\sqrt{\frac{1-v}{1+v}}}-\frac{u}{1-v}S\parens{u\sqrt{\frac{1+v}{1-v}}}}du.
\end{equation}
\normalsize

We now focus on the term in square brackets in the above expression. We write its power series representation as:

\begin{equation} \label{niceone5}
\bracks{\frac{u}{1+v}S\parens{u\sqrt{\frac{1-v}{1+v}}}-\frac{u}{1-v}S\parens{u\sqrt{\frac{1+v}{1-v}}}}= \sum_{r=1}^\infty s_r(u)v^{2r-1},
\end{equation}
for a sequence of functions $s_r(u), r\geq 1$, which are not easy to obtain. Note that the right hand side of Equation (\ref{niceone5}) is analytic in $v$ near $v=0$ and is also odd in $v$. Placing the right hand side of Equation (\ref{niceone5}) into expression (\ref{niceone2}), and interchanging integrals, expression (\ref{niceone2}) becomes:

\begin{equation} \label{almostthere}
\sum_{r=1}^\infty \int_0^\infty 2S(u)T(u)s_r(u)du\int_{z}^1\frac{v}{\sqrt{v^2-z^2}}v^{2(r-1)}dv.
\end{equation}
Substituting $v^2=t^2+z^2$, and noting that $vdv=tdt$, the second integral in expression (\ref{almostthere}) becomes:  
\beqn
\int_0^{\sqrt{1-z^2}}(z^2+t^2)^{r-1}dt= \sqrt{1-z^2}\int_0^1\parens{z^2(1-v^2)+v^2}^{r-1}dv,
\eeqn
where we have substituted $t=v\sqrt{1-z^2}$. This is quite convenient since the factor $\sqrt{1-z^2}$ cancels on both sides of Equation (\ref{massive2}). The final result is the slightly simpler identity below:

\footnotesize
\begin{equation} \label{identity}
\sum_{n=1}^\infty \frac{z^{2(n-1)}}{2n} \expe{\theta^{2n}}\frac{(n!)^22^{2n}}{(2n)!}= \sum_{r=1}^\infty \sum_{k=0}^{r-1} z^{2k} \int_0^1 {r-1\choose k} (1-v^2)^kv^{2(r-1-k)}dv\int_0^\infty2S(u)T(u)s_r(u)du.
\end{equation}
\normalsize
Comparing coefficients of powers of $z$ \textit{gives the moments} of the distribution of $\theta$. Collecting terms with  $k=n-1$, we can simplify Equation (\ref{identity}) to obtain the following theorem:

\begin{theorem} \label{thm3}
\small
\begin{equation} \label{niceone6}
\expe{\theta^{2n}}= {2n \choose n} \frac{2n}{2^{2n}}\sum_{r=n}^\infty \int_0^1 {r-1\choose n-1} (1-v^2)^{n-1}v^{2(r-n)}dv\int_0^\infty2S(u)T(u)s_r(u)du,
\end{equation}
where the $s_r$ are implicitly determined by Equation (\ref{niceone5}).
\end{theorem}

\normalsize
Proposition \ref{prop10} below gives an even more explicit formula for the second moment of $\theta$:
\begin{proposition} \label{prop10}
The second moment of $\theta$, corresponding to $n=1$, can be calculated explicitly as 
\footnotesize
\beqn
\int_0^{\infty} du_1\int_0^{u_1}\sqrt{\frac{u_1u_2}{\mathrm{sinh}\, u_1 \mathrm{sinh}\, u_2}}\frac{u_1u_2}{u_1+u_2} \parens{\frac{\frac{1}{u_1^2}\parens{1-\frac{u_1}{\mathrm{sinh} u_1}\mathrm{cosh} u_1}-\frac{1}{u_2^2}\parens{1-\frac{u_2}{\mathrm{sinh}\, u_2} \mathrm{cosh}u_2}}{u_1-u_2}} du_2 .
\eeqn
\end{proposition}
\normalsize
\begin{proof}
For the case $n=1$, we simplify Equation (\ref{niceone6}) as:

\begin{equation} \label{theone}
\expe{\theta^{2}}=2\int_0^\infty S(u)T(u)du\int_0^1 v^{-1}dv\sum_{r=1}^\infty v^{2r-1}s_r(u).
\end{equation}
Placing the left hand side of Equation (\ref{niceone5}) into the above equality, we have $\expe{\theta^{2}}$ equals

\footnotesize
\begin{equation} \label{cooldude}
\begin{split}
2\int_0^\infty du\int_0^1v^{-1}S(u)T(u)\bracks{\frac{u}{1+v}S\parens{u\sqrt{\frac{1-v}{1+v}}}-\frac{u}{1-v}S\parens{u\sqrt{\frac{1+v}{1-v}}}}dv.
\end{split}
\end{equation}
\normalsize
Making the transformation $u \rightarrow u\sqrt{1+v}$ in the first term of the bracketed expression above and $u \rightarrow u\sqrt{1-v}$ in the second term of the bracketed expression above, expression (\ref{cooldude}) becomes

\footnotesize

\beqn
2\int_0^\infty udu\int_0^1v^{-1}\,S(u\sqrt{1+v})S(u\sqrt{1-v})\parens{T(u\sqrt{1+v})-T(u\sqrt{1-v})}dv.
\eeqn
\normalsize
We now make one final transformation. Let $u_1=u\sqrt{1+v}$ and let $u_2=u\sqrt{1-v}$. The Jacobian of this transformation is $\frac{u}{\sqrt{1-v^2}}$. Note that $v=\frac{u_1^2-u_2^2}{u_1^2+u_2^2}$, and $\sqrt{1-v^2}=\frac{u_1u_2}{u^2}$. After making this transformation, we obtain a simpler form of expression (\ref{cooldude}), which gives us an explicit formula for the second moment of $\theta$:

\footnotesize

\begin{equation} \label{finaldude}
\expe{\theta^2}=\int_0^{\infty} du_1 \int_0^{u_1}\parens{\frac{2u_1u_2}{u_1+u_2}}S(u_1)S(u_2)\frac{T(u_1)-T(u_2)}{u_1-u_2}du_2.
\end{equation}
\normalsize
Using the definitions of $S$ and $T$, we arrive at our final expression for the second moment of $\theta$ as:

\footnotesize
\beqn
\int_0^{\infty} du_1\int_0^{u_1}\sqrt{\frac{u_1u_2}{\text{sinh}\, u_1 \text{sinh}\, u_2}}\frac{u_1u_2}{u_1+u_2} \parens{\frac{\frac{1}{u_1^2}\parens{1-\frac{u_1}{\text{sinh} u_1}\text{cosh}\, u_1}-\frac{1}{u_2^2}\parens{1-\frac{u_2}{\text{sinh}\, u_2} \text{cosh}\,u_2}}{u_1-u_2}} du_2 .
\eeqn
\normalsize
This finishes the proof.
\end{proof}
Using the above expression for the second moment of $\theta$ we numerically obtain a value of .240522. Recall that the Monte Carlo simulation in Table \ref{tab:conditional_prob_sim} reported a value of .235057.
\begin{remark}
One must proceed numerically using Equation (\ref{niceone6}) to calculate higher order moments. 
\end{remark}

\subsection{Open Problems}
This work admits many potential extensions. Other Gaussian processes should be treatable with our methods, for example, a pair of Ornstein-Uhlenbeck processes or a pair of \textit{correlated} Wiener processes. It would be of particular interest to compare the variances of the correlation coefficient amongst these cases.

\section{Acknowledgments} First and foremost, we are very grateful for the very careful and detailed work of the anonymous referees, whose detailed and very careful reports significantly improved the quality of this work. We thank Professor Amir Dembo for his advice. We are very grateful to Professor Lawrence Brown for his invaluable input. We are also extremely grateful to Professor Edward George for his careful reading of this manuscript, for his invaluable suggestions, and, most importantly, for suggesting ``volatile correlation'' to describe the phenomenon of ``nonsense'' correlation.

\newpage

\bibliographystyle{plain}
\bibliography{refs2}

\begin{thebibliography}{10}

\bibitem{Aldrich}
J.~Aldrich.
\newblock Correlations genuine and spurious in pearson and yule.
\newblock {\em Stat Sci}, 10:\,364--376, 1995.

\bibitem{boas}
R.~Boas.
\newblock {\em Entire Functions}.
\newblock Academic Press, New York, 1954.

\bibitem{Kac}
P.~Erd{\"o}s and M.~Kac.
\newblock On certain limit theorems of the theory of probability.
\newblock {\em Bull. Amer. Math. Soc.}, 52:\,292--302, 1946.

\bibitem{feller}
W.~Feller.
\newblock {\em An Introduction to Probability Theory and Its Applications}.
\newblock John Wiley, New York, 1953.

\bibitem{gran}
C~Granger and D.~Newbold.
\newblock Spurious regression in econometrics.
\newblock {\em J. Econometrics}, 2:\,111--120, 1974.

\bibitem{Hendry}
D.F. Hendry.
\newblock Economic modelling with cointegrated variables: an overview.
\newblock {\em Oxf Bull Econ Stat}, 48:\,201--212, 1986.

\bibitem{shepp1}
B~Logan, C~Mallows, S~Rice, and L~Shepp.
\newblock Limit distributions of self-normalized sums.
\newblock {\em Ann. Probab.}, 1:788--809, 1973.

\bibitem{shepp2}
B~Logan and L~Shepp.
\newblock Real zeros of random polynomials.
\newblock {\em Proc. London Math. Soc.}, 18:\,29--35, 1968.

\bibitem{magnus}
J.~Magnus.
\newblock The exact moments of a ratio of quadratic forms in normal variables.
\newblock {\em Ann. Econ. Statist.}, pages \,95--109, 1986.

\bibitem{Mann}
M.~Mann, Z.~Zhang, M.~K. Hughes, R.~S. Bradley, S.~K. Miller, S.~Rutherford,
  and F.~Ni.
\newblock Proxy-based reconstructions of hemispheric and global surface
  temperature variations over the past two millenia.
\newblock {\em Proc. Natl. Acad. Sci.}, 105:13252--13257, 2008.

\bibitem{McShane}
Blakeley~B. McShane and Abraham~J. Wyner.
\newblock A statistical analysis of multiple temperature proxies: Are
  reconstructions of surface temperatures over the last 1000 years reliable?
\newblock {\em Ann. Appl. Stat.}, 5(1):5--44, 03 2011.

\bibitem{phil}
P.C.B. Phillips.
\newblock Understanding spurious regressions in econometrics.
\newblock {\em J. Econometrics}, 33:\,311--340, 1986.

\bibitem{shepp3}
L.C.G. Rogers and L~Shepp.
\newblock The correlation of the maxima of correlated {B}rownian motion.
\newblock {\em J. Appl. Probab.}, 43:\,880--883, 2006.

\bibitem{shepp4}
L.~Shepp.
\newblock The joint density of the maximum and its location for a {W}iener
  process with drift.
\newblock {\em J. Appl. Probab.}, 16:\,423--427, 1979.

\bibitem{Yule}
G.~U. Yule.
\newblock Why do we sometimes get nonsense correlations between time series? a
  study in sampling and the nature of time series.
\newblock {\em Roy. Statist. Soc.}, 89:\,1--63, 1926.

\end{thebibliography}

\end{document}